\documentclass[10pt]{article}
\usepackage{epsfig,amssymb,amsmath}
\usepackage{palatino}
\usepackage{verbatim}
\usepackage{enumerate}

\def\baa {\begin{eqnarray*}}
\def\eaa {\end{eqnarray*}}

\def \la {\lambda}

\def\la{\lambda}

\def \sign {\rm sign\,}
\newtheorem{lemma}{Lemma}[section]

\newtheorem{corollary}[lemma]{Corollary}
\newtheorem{theorem}[lemma]{Theorem}

\def\bc  {}
\makeatletter \@addtoreset{equation}{section} \makeatother
\begin{document}

\title{Tur\'{a}n's inequality for ultraspherical polynomials revisited}

\author{Geno Nikolov}

\date{}
\maketitle

\begin{abstract}
We present a short proof that the normalized Tur\'{a}n determinant
in the ultraspherical case is convex or concave depending on whether
parameter $\lambda$ is positive or negative.
\end{abstract}

\smallskip



\section{Introduction and statement of the result}

In the 40's of the last century, while studying the zeros of
Legendre polynomials $P_n(x)$, P. Tur\'{a}n discovered the
inequality
\begin{equation}\label{e1.1}
P_n^{2}(x)-P_{n-1}(x)P_{n+1}(x)\ge 0, \ \ -1\le x\le 1,
\end{equation}
with equality only for $x=\pm 1$. Since the left-hand side of
\eqref{e1.1} is representable in determinant form,
$$
\Delta_n(x)=\begin{vmatrix} P_n(x) & P_{n+1}(x)\\
P_{n-1}(x) & P_n(x) \end{vmatrix}
$$
$\Delta_n(x)$ is referred to as \emph{Tur\'{a}n's determinant}.

The result of Tur\'{a}n inspired a considerable interest, and by now
there is a vast amount of publications on the so-called
\emph{Tur\'{a}n type inequalities}. G. Szeg\H{o} \cite{GS1948} gave
four different proof of \eqref{e1.1}. As Szeg\H{o} pointed out in
\cite{GS1948}, his third proof extends Tur\'{an}'s inequality to
other classes of functions including ultraspherical polynomials,
Laguerre and Hermite polynomials, Bessel functions, etc. This idea
was rediscovered by Skovgaard \cite{HS1954}.

Karlin and Szeg\H{o} \cite{KS1961} posed the problem of
characterizing the set of pairs $\{\alpha,\beta\}$ for which the
normalized Jacobi polynomials
$P_m^{(\alpha,\beta)}(x)/P_m^{(\alpha,\beta)}(1)$ admit a Tur\'{a}n
type inequality. Szeg\H{o} proved that Tur\'{a}n's inequality holds
whenever $\beta\geq|\alpha|$, $\alpha>-1$. In two subsequent papers
G. Gasper \cite{GG1971, GG1972} improved Szeg\H{o}'s result showing
finally that the sought pairs $\{\alpha,\beta\}$ are those
satisfying $\beta\ge\alpha>-1$.

Our concern here is Tur\'{an}'s inequality in the ultraspherical
case. Throughout this paper, $p_n^{(\la)}$ stands for the $n$-th
ultraspherical polynomial normalized to assume value $1$ at $x=1$,
$$
p_n^{(\la)}(x)=\frac{P_n^{(\la)}(x)}{P_n^{(\la)}(1)}\,.
$$
Let
\begin{equation}\label{e1.2}
\Delta_{n,\la}(x):=\big[p_n^{(\la)}(x)\big]^2-p_{n-1}^{(\la)}(x)
p_{n+1}^{(\la)}(x)\,,
\end{equation}
then Tur\'{a}n's inequality for ultraspherical polynomials reads as
\begin{equation}\label{e1.3}
\Delta_{n,\la}(x)\geq 0,\qquad x\in [-1,1]\,.
\end{equation}
To the many proofs of \eqref{e1.3} (see, e.g. \cite{AD1955, OS1951,
GS1948, TN1951, VLR1957}), let us add the one in \cite{GN2013} based
on a Hermite interpolation formula, yielding the representation
$$
\Delta_{n,\lambda}(x)=\frac{1-x^2}{n(n+2\lambda)}\,\sum_{k=1}^{n}\ell_k^{2}(x)
(1-x_kx)\, \big[p_n^{\prime}(x_k)\big]^2
$$
(here, $\{\ell_k\}_{k=1}^{n}$ are the Lagrange basis polynomials for
interpolation at the zeros $\{x_k\}_{k=1}^{n}$ of
$p_n=p_n^{(\la)}$).

Since $\Delta_{n,\la}(\pm 1)=0$, it is of interest to describe the
behavior of the normalized Tur\'{a}n function
\begin{equation}\label{e1.4}
\varphi_{n,\la}(x):=\frac{\Delta_{n,\la}(x)}{1-x^2}\,.
\end{equation}
Thiruvenkatachar and Nanjundiah \cite{TN1951} have shown that
$\varphi_{n,\la}$ increases in $[-1,0]$ and decreases in $[0,1]$ if
$-1/2<\la<0$, and has the opposite behavior if $\la>0$. Since
$\varphi_{n,\la}$ is an even function, it follows that for $x\in
[-1,1]$,
\begin{equation*}
\begin{split}
&\varphi_{n,\la}(1)\leq\varphi_{n,\la}(x)\leq \varphi_{n,\la}(0),
\quad -1/2<\la<0 \\
&\varphi_{n,\la}(0)\leq\varphi_{n,\la}(x)\leq \varphi_{n,\la}(1),
\quad \la>0,
\end{split}
\end{equation*}
which together with $\varphi_{n,\la}(0)=\Delta_{n,\la}(0)$ and
$\varphi_{n,\la}(1)=-\Delta_{n,\la}^{\prime}(1)/2=1/(2\la+1)$ yields
the following two-sided estimates for $\Delta_{n,\la}(x)$ when $x\in
[-1,1]$.
\begin{equation}\label{e1.5}
\begin{split}
&\frac{1-x^2}{2\la+1}\leq\Delta_{n,\la}(x)\leq
\Delta_{n,\la}(0)(1-x^2),
\quad -1/2<\la<0 \\
&\Delta_{n,\la}(0)(1-x^2)\leq\Delta_{n,\la}(x)\leq
\frac{1-x^2}{2\la+1}, \quad \la>0.
\end{split}
\end{equation}
Here we make this observation more precise by proving the following:

\begin{theorem}\label{t1.1}
The normalized Tur\'{a}n function $\varphi_{n,\la}$ is concave or
convex on $\mathbb{R}$ depending on whether $-1/2<\la<0$ or $\la>0$.
\end{theorem}
(Note that  $\varphi_{n,0}\equiv 1$.) Theorem~\ref{t1.1} implies
two-sided estimates for the Tur\'{a}n determinant
$\Delta_{n,\la}(x)$ which both sharpen and extend \eqref{e1.5} for
all $x\in \mathbb{R}$.
\begin{corollary}\label{c1.2}
\begin{itemize}
\item[(i)] If $-1/2<\la<0$, then
$$
\Big[(1-|x|)\Delta_{n,\la}(0)+\frac{|x|}{2\la+1}\Big](1-x^2)
\leq\Delta_{n,\la}(x)\leq \Delta_{n,\la}(0)(1-x^2)\,,\quad x\in
\mathbb{R}.
$$
\item[(ii)] If $\la>0$, then
$$
\Delta_{n,\la}(0)\,(1-x^2)\leq \Delta_{n,\la}(x)\leq
\Big[(1-|x|)\Delta_{n,\la}(0)+\frac{|x|}{2\la+1}\Big]\,(1-x^2)\,,
\quad x\in \mathbb{R}.
$$
\end{itemize}
\end{corollary}

The proof of Theorem~\ref{t1.1} is given in the next section. The
last section contains some remarks and comments.

\section{Proof}
We shall work with the renormalized ultraspherical polynomials
$$
p_n^{(\la)}(x)=P_n^{(\la)}(x)/P_n^{(\la)}(x)\,,
$$
and for simplicity's sake we omit the superscript $^{(\la)}$, so
$p_n:=p_n^{(\la)}$. The next two identities readily follow from
\cite[equation (4.7.28)]{Sze75}:
\begin{eqnarray*}
&&p_n(x)=-\frac{1}{n+2\la}\,x\,p_n^{\prime}(x)
+\frac{1}{n+1}\,p_{n+1}^{\prime}(x)\,,\\
&&p_{n+1}(x)=-\frac{1}{n+2\la}\,p_n^{\prime}(x)
+\frac{1}{n+1}\,x\,p_{n+1}^{\prime}(x)\,.
\end{eqnarray*}
These identities are used for deriving representations of $p_{n+1}$
and $p_{n-1}$ in terms of $p_n$ and $p_n^{\prime}$:
\begin{eqnarray*}
&&p_{n+1}(x)=x\,p_n(x)-\frac{1-x^2}{n+2\la}\,p_n^{\prime}(x)\,,\\
&&p_{n-1}(x)=x\,p_n(x)+\frac{1-x^2}{n}\,p_n^{\prime}(x)\,.
\end{eqnarray*}
By replacing $p_{n+1}$ and $p_{n-1}$ in
$\Delta_{n,\la}=p_n^2-p_{n-1}p_{n+1}$ we obtain
$$
\Delta_{n,\la}(x)=\frac{1-x^2}{n(n+2\la)}\,\Big[n(n+2\la)p_n^2(x)
-2\la\,x\,p_n(x)p_n^{\prime}(x)+(1-x^2)\big[p_n^{\prime}(x)\big]^2\Big]\,,
$$
hence
\begin{equation}\label{e2.1}
\varphi_{n,\la}(x)=\frac{1}{n(n+2\la)}\,\Big[n(n+2\la)p_n^2(x)
-2\la\,x\,p_n(x)p_n^{\prime}(x)+(1-x^2)\big[p_n^{\prime}(x)\big]^2\Big]\,.
\end{equation}
Differentiating \eqref{e2.1} and using  the differential equation
\begin{equation}\label{e2.2}
(1-x^2)y^{\prime\prime}-(2\la+1)x\,y^{\prime}+n(n+2\la)\,y=0\,,
\qquad y=p_n(x),
\end{equation}
we find
\begin{equation}\label{e2.3}
\begin{split}
\varphi_{n,\la}^{\prime}(x)&=\frac{2\la}{n(n+2\la)}\Big[ x
\big[p_n^{\prime}(x)\big]^2-p_n(x)p_n^{\prime}(x)
-x\,p_n(x)p_n^{\prime\prime}(x)\Big]\\
&=-\frac{2\la}{n(n+2\la)}\,p_n^2(x)\,
\Big(\frac{x\,p_n^{\prime}(x)}{p_n(x)}\Big)^{\prime}\,.
\end{split}
\end{equation}
Let $x_1<x_2<\cdots<x_n$ be the zeros of $p_n$, they form a
symmetric set with respect to the origin, therefore
$$
\frac{p_n^{\prime}(x)}{p_n(x)}=\sum_{k=1}^{n}\frac{1}{x-x_k}
=\frac{1}{2}\,\sum_{k=1}^{n}\Big(\frac{1}{x-x_k}+\frac{1}{x+x_k}\Big)
=x\,\sum_{k=1}^{n}\frac{1}{x^2-x_k^2}\,.
$$
Consequently,
$$
\Big(\frac{x\,p_n^{\prime}(x)}{p_n(x)}\Big)^{\prime}=
-2x\,\sum_{k=1}^{n}\frac{x_k^2}{(x^2-x_k^2)^2}\,,
$$
and \eqref{e2.3} implies
\begin{equation}\label{e2.4}
\varphi_{n,\la}^{\prime}(x)=\frac{4\la\,x}{n(n+2\la)}
\,\sum_{k=1}^{n}x_k^2\,q_{n,k}^2(x)\,,\qquad
q_{n,k}(x)=\frac{p_n(x)}{x^2-x_k^2}\,.
\end{equation}
Now \eqref{e2.4} shows that $\sign \varphi_{n,\la}^{\prime}(x)=
\sign \la\,x$, a result already obtained by Thiruvenkatachar and
Nanjundiah \cite{TN1951}. In fact, \eqref{e2.4} implies more than
that, namely,
\begin{equation}\label{e2.5}
\sign \varphi_{n,\la}^{(r)}(x)=\sign \la\,,\qquad x>x_n,\ \
r=1,2,\ldots, 2n-2\,.
\end{equation}
Indeed, $\varphi_{n,\la}^{\prime}$ is a sum of polynomials with
leading coefficients of the same sign as $\la$ and with all their
zeros being real and located in $[x_1,x_n]$. By Rolle's theorem, the
derivatives of these polynomials inherit the same properties, hence
they have no zeros in $(x_n,\infty)$ and therefore have the same
sign as $\la$ therein. In particular, \eqref{e2.5} implies
$$
\sign \varphi_{n,\la}^{\prime\prime}(x)=\sign \la,\qquad x\in
(x_n,\infty)
$$
and to prove Theorem~\ref{t1.1} we need to show that $\sign
\varphi_{n,\la}^{\prime\prime}(x)=\sign \la$ for $x\in (0,x_n]$. In
view of \eqref{e2.3}, this is equivalent to prove that the function
\begin{equation}\label{e2.6}
\psi_{n,\la}(x):=\big[p_n^{\prime}(x)\big]^2-p_n(x)p_n^{\prime}(x)
-x\,p_n(x)p_n^{\prime\prime}(x)
\end{equation}
satisfies
\begin{equation}\label{e2.7}
\psi_{n,\la}^{\prime}(x)>0,\qquad x\in (0,x_n].
\end{equation}
We differentiate \eqref{e2.6} and make use of the differential
equations \eqref{e2.2} and
$$
(1-x^2)y^{\prime\prime\prime}-(2\la+3)x\,y^{\prime\prime}
+(n-1)(n+2\la+1)\,y^{\prime}=0\,, \qquad y=p_n(x),
$$
to obtain a representation of $\psi_{n,\la}^{\prime}(x)$ as a
quadratic form of $p_n^{\prime}$ and $p_n^{\prime\prime}$:
\begin{equation*}
\begin{split}
n(n+2\la)(1-x^2)\psi_{n,\la}^{\prime}(x)=&(2\la+1)(n-1)(n+2\la+1)x^2\,
\big[p_n^{\prime}(x)\big]^2\\
&-(2\la+1)x\big[1+2(\la+1)x^2\big]\,p_n^{\prime}(x)p_n^{\prime\prime}(x)\\
&+(1-x^2)\big[2+(2\la+1)x^2\big]\,\big[p_n^{\prime\prime}(x)\big]^2
\end{split}
\end{equation*}
The discriminant $D$ of this quadratic form equals
$$
D(x)=(2\la+1)\,x^2\big[2\la+3-(2\la+1)(1-x^2)\big]\,D_1(x),
$$
where
$$
D_1(x)=(2\la+1)\,\frac{\big[2\la+3-(2\la+2)(1-x^2)\big]^2}{2\la+3-(2\la+1)(1-x^2)}
-4(n-1)(n+2\la+1)(1-x^2)\,.
$$
Our goal is to prove that
\begin{equation}\label{e2.8}
D_1(x)<0\,,\qquad x\in (0,x_n],
\end{equation}
which implies $D(x)<0$ and consequently $\psi_{n,\la}^{\prime}(x)>0$
in $(0,x_n]$. It is readily verified that
$$
\frac{\big[2\la+3-(2\la+2)(1-x^2)\big]^2}{2\la+3-(2\la+1)(1-x^2)}
\leq 2\la+3-(2\la+5/2)(1-x^2)\,, \qquad x\in [-1,1],
$$
therefore
\begin{equation*}
\begin{split}
D_1(x)&\leq
(2\la+1)\,\big[2\la+3-(2\la+5/2)(1-x^2)\big]-4(n-1)(n+2\la+1)(1-x^2)\\
&=(2\la+1)(2\la+3)-\big[4(n+\la)^2-(\la+3/2)\big](1-x^2)\,, \qquad
x\in [-1,1].
\end{split}
\end{equation*}
Hence, to prove \eqref{e2.8}, it suffices to show that
$$
1-x^2>\frac{(2\la+1)(2\la+3)}{4(n+\la)^2-\la-3/2}\,,\qquad x\in
(0,x_n)
$$
or, equivalently,
\begin{equation}\label{e2.9}
x_n^2<1-\frac{(2\la+1)(2\la+3)}{4(n+\la)^2-\la-3/2}\,.
\end{equation}
Thus, we need an upper bound for $x_n$, the largest zero of the
ultraspherical polynomial $P_n^{(\la)}$. Amongst the numerous upper
bounds in the literature, we use the one from \cite[Lemma~6]{GN2005}
(see also \cite[p. 1801]{DN2010}):
\begin{equation}\label{e2.10}
x_n^2<\frac{(n+\la)^2-(\la+1)^2}{(n+\la)^2+3\la+5/4+3(\la+1/2)^2/(n-1)}\,.
\end{equation}
The comparison of the right-hand sides of \eqref{e2.9} and
\eqref{e2.10} (we have used \emph{Wolfram Mathematica} for this
purpose) shows that the latter is the smaller one, hence
\eqref{e2.9} holds true. With this \eqref{e2.7} is proved, hence
$\sign \varphi_{n,\la}^{\prime\prime}(x)=\sign \la$ for $x\in
(0,x_n]$ and consequently
$$
\sign \varphi_{n,\la}^{\prime\prime}(x)=\sign \la,\qquad x\in
(0,\infty).
$$
Since $\varphi_{n,\la}^{\prime\prime}$ is an even function, this
accomplishes the proof of Theorem~\ref{t1.1}.

\section{Remarks}

\noindent 1) There are also some results concerning concavity of
$\Delta_{n,\la}$. In the classical Tur\'{a}n case, $\la=1/2$,
Madhava Rao and Thiruvenkatachar \cite{RT1949} proved that
$$
\frac{d^2}{dx^2}\,\Delta_n(x)
=-\frac{2}{n(n+1}\,\big[P_n^{\prime\prime}(x)]^2\,,
$$
showing that $\Delta_n$ is a concave function. Venkatachaliengar and
Lakshmana Rao \cite{VLR1957} extended this result by proving that
$\Delta_{n,\la}$ is a concave function in $[-1,1]$ provided  $\la\in
(0,1/2]$. Generally, $\Delta_{n,\la}$ is neither convex nor concave
if $\la\not\in [0,1/2]$.\medskip

\noindent 2). Sz{\'a}sz \cite{OS1951} proved the following pair of
bounds for $\Delta_{n,\la}(x)$:
$$
\frac{\la\big(1-[p_n^{(\la)}(x)]^2\big)}{(n+\la-1)(n+2\la)}<
\Delta_{n,\la}(x) <
\frac{n+\la}{\la+1}\,\frac{\Gamma(n)\Gamma(2\la+1)}{\Gamma(n+2\la+1)}\,,
\qquad \la\in (0,1)\,.
$$
\smallskip

\noindent 3). In a recent paper \cite{NP2015} we gave both an
analytical and a computer proof of the following refinement of
Tur\'{a}n's inequality:
$$
|x|\,\big[p_n^{(\la)}(x)\big]^2-p_{n-1}^{(\la)}(x)p_{n+1}^{(\la)}(x)\geq
0, \qquad x\in [-1,1],\ -1/2<\la\leq 1/2\,,
$$
with the equality occurring only for $x=\pm 1$ and, if $n$ is even,
$x=0$. This inequality provides another lower bound for
$\Delta_{n,\la}(x)$ in the case $-1/2<\la\leq 1/2$. A computer proof
of the Legendre case ($\la=1/2$) was given earlier by Gerhold and
Kauers \cite{GK2006}.
\medskip

\noindent 4). In \cite{TN1951} the authors proved also monotonicity
of $\Delta_{n,\la}(x)$, $x\in [-1,1]$ fixed, with respect to $n$. We
refer to \cite{BS2007, RS1995} for some general condition on the
sequences defining the three-term recurrence relation for orthogonal
polynomials, which ensure the monotonicity of the associated
Tur\'{a}n determinants.
\medskip

\noindent 5). For a higher order Tur\'{an} inequalities and a
discussion on the interlink between the Tur\'{a}n type inequalities
and the Riemann hypothesis or the recovery of the orthogonality
measure, we refer to \cite{DD1998} and the references therein.

\bibliographystyle{amsplain}

\bigskip

\noindent
{\sc Geno Nikolov} \smallskip\\
Department of Mathematics and Informatics\\
Universlty of Sofia \\
5 James Bourchier Blvd. \\
1164 Sofia \\
BULGARIA \\
{\it E-mail:} {\tt geno@fmi.uni-sofia.bg}

\end{document}